\documentclass[12pt,reqno]{amsart}
\usepackage{a4,latexsym,amssymb,amsfonts,pgf}
\usepackage{enumerate}
\usepackage{amsthm}
\usepackage{amsmath}
\usepackage[francais,english]{babel}
\usepackage{hyperref}
\usepackage[active]{srcltx}
\usepackage[latin1]{inputenc}
\usepackage{amscd}
\usepackage{pstricks}
\usepackage{dsfont}
\usepackage[all]{xy}
\usepackage{float}
\usepackage{geometry}
\newtheorem{thm}{Th\'eor\`eme}
\newtheorem*{thmp}{Théorème principal}
\newtheorem{Def}{D\'efinition}
\newtheorem{prop}{Proposition}
\newtheorem{cor}{Corollaire}

\newtheorem{exs}{Exemples}

\newtheorem{lem}{Lemme}
\newtheorem{rem}{Remarque}
\newtheorem{rems}{Remarques}
\newcommand{\bexs}{\begin{exs}}
\newcommand{\eexs}{\end{exs}}
\newcommand{\bt}{\begin{thm}}
\newcommand{\et}{\end{thm}}
\newcommand{\bd}{\begin{Def}}
\newcommand{\ed}{\end{Def}}
\newcommand{\bl}{\begin{lem}}
\newcommand{\el}{\end{lem}}
\newcommand{\bp}{\begin{prop}}
\newcommand{\ep}{\end{prop}}
\newcommand{\brem}{\begin{rem}}
\newcommand{\erem}{\end{rem}}
\newcommand{\brems}{\begin{rems}}
\newcommand{\erems}{\end{rems}}
\newcommand{\bc}{\begin{cor}}
\newcommand{\ec}{\end{cor}}
\newcommand{\bpr}{\begin{proof}}
\newcommand{\epr}{\end{proof}}
\def\NN{\mathbb{N}}
\def\RR{\mathbb{R}}

\def\QQ{\mathbb{Q}}

\def\ZZ{\mathbb{Z}}
\begin{document}
\title[ Capitulation des 2-classes d'idéaux ]{ Capitulation
des 2-classes d'idéaux de $\mathbf{k}=\QQ(\sqrt{ 2p}, i)$}


\author[Abdelmalek AZIZI]{Abdelmalek AZIZI$^\dag$}

\thanks{$^\dag$Recherche soutenue par l'Académie Hassan II des Sciences et
Techniques, Maroc}

\address{Département de Mathématiques, Faculté des Sciences,Université Mohammed 1, Oujda, Maroc }

\author{Mohammed TAOUS}


\email{abdelmalekazizi@yahoo.fr}
\email{taousm@hotmail.com}



\maketitle

\begin{abstract}
Let $p$ be a prime number such that $p\equiv 1$ mod $8$ and
$i=\sqrt{-1}$. Let $\mathbf{k}=\QQ(\sqrt{2p}, i)$,
$\mathbf{k}_1^{(2)}$ be the Hilbert $2$-class field of
$\mathbf{k}$, $\mathbf{k}_2^{(2)}$ be the Hilbert $2$-class field
of $\mathbf{k}_1^{(2)}$ and $G=\mathrm{Gal}(\mathbf{k}_2^{(2)}/\mathbf{k})$
be the Galois group of $\mathbf{k}_2^{(2)}/\mathbf{k}$. Suppose
that the $2$-part, $C_{\mathbf{k}, 2}$, of the class group of
$\mathbf{k}$ is of type $(2, 4)$; then $\mathbf{k}_1^{(2)}$
contains six extensions $\mathbf{K_{i,  j}}/\mathbf{k}$, $i=1, 2,
3$ and $j=2, 4$. Our goal is to study the problem of the
capitulation of $2$-ideal classes of $\mathbf{K_{i,  j}}$ and to
determine the structure of $G$.
\end{abstract}

\selectlanguage{francais}

\begin{abstract}
Soient $p$ un nombre premier tel que $p\equiv 1$ mod $8$ et
$i=\sqrt{-1}$. Soient $\mathbf{k}=\QQ(\sqrt{2p}, i)$,
$\mathbf{k}_1^{(2)}$ le $2$-corps de classes de Hilbert de
$\mathbf{k}$, $\mathbf{k}_2^{(2)}$ le $2$-corps de classes de
Hilbert de $\mathbf{k}_1^{(2)}$ et
$G=Gal(\mathbf{k}_2^{(2)}/\mathbf{k})$ le groupe de Galois de
$\mathbf{k}_2^{(2)}/\mathbf{k}$. Supposons que la $2$-partie,
$C_{\mathbf{k}, 2}$, du groupe de classes de $\mathbf{k}$ est de
type $(2, 4)$; alors $\mathbf{k}_1^{(2)}$ contient six extensions
$\mathbf{K_{i,  j}}/\mathbf{k}$, $i=1, 2, 3$ et $j=2, 4$. On
s'intéresse au problème de capitulation des $2$-classes de
$\mathbf{k}$ dans $\mathbf{K_{i, j}}$ et à déterminer la structure
de $G$.
\end{abstract}
\section{Introduction}
Soient $k$ un corps de nombres de degré fini sur $\QQ$, $p$ un nombre
premier, $C_k$ le groupe de classes de $k$ et $C_{k, p}$ le
$p$-groupe de classes de $k$. On note $k_p^{(1)}$ le $p$-corps de
classes de Hilbert de $k$ au sens large. Soit  $k_p^{(n)}$ (pour $n$
un entier naturel) la suite de $p$-corps de classes de Hilbert définie
par: $k_p^{(0)}=k$ et $k_p^{(n+1)}=(k_p^{(n)})_p^{(1)}$. Alors on a
$$k_p^{(0)}\subseteq k_p^{(1)}\subseteq ....\subseteq k_p^{(n)}\subseteq ...$$
\indent Cette suite est appelée la tour des $p$-corps de classes de
Hilbert de $k$, on sait qu'elle est finie si et seulement s'il
existe une $p$-extension finie $E$ de $k$ telle que le $p$-nombre de
classes (la $p$-partie du nombre de classes) de $E$ est égale à $1$.
Mais cette caractérisation ne permet pas d'obtenir une procédure
pour dire que cette suite s'arrête ou non, cependant, il est connu
par un résultat de groupes de Taussky (\cite{Ta37}) que si
$C_{k_{p}^{(1)}, p}$ est cyclique alors $C_{k_{p}^{(2)}, p}$ est
trivial, ce qui implique que $k_p^{(2)}=k_p^{(3)}$. Or, si $p=2$ et
$C_{k, p}$ est de type $(2, 4)$, alors d'après un autre résultat de
groupes de Blackburn
(\cite{Bl58}), on a que le rang de $C_{k_{p}^{(1)}, p}\leq 3$.\\
\indent L'objet de ce travail est l'étude du problème de la tour
pour le corps $\mathbf{k}=\QQ(\sqrt{2p}, i)$, dont le $2$-groupe
de classes est de type $(2, 4)$. Nous déterminons aussi les
$2$-classes de $C_{\mathbf{k}}$ qui capitulent dans les
sous-extensions propres de $\mathbf{k}_2^{(1)}/\mathbf{k}$, ce qui
nous permet de trouver une représentation de
$G=Gal(\mathbf{k}_2^{(2)}/\mathbf{k})$ groupe de Galois de
$\mathbf{k}^{(2)}_2/\mathbf{k}$, lorsque $\mathbf{k}_2^{(1)}\neq
\mathbf{k}_2^{(2)}$. Notre théorème principal est le suivant :
\begin{thmp}
Soit $\mathbf{k}=\QQ(\sqrt{2p}, i)$ avec $p$ un nombre premier tel
que $p\equiv 1\mod 8$ et $(\frac{2}{p})_4=(\frac{p}{2})_4=-1$.
Alors il existe deux entiers $e$, $f\in\NN$ tels que
$p=e^2+16f^2$. Soient $\pi_1=e+4fi$, $\pi_2=e-4fi$, $2^n$ le
$2$-nombre de classes de $\QQ(\sqrt{-p})$, $\mathcal{H}_1$,
$\mathcal{H}_2$ et $\mathcal{H}$ les idéaux premiers au-dessus de
$\pi_1$, $\pi_2$ et $1+i$ dans $\mathbf{k}$. Alors
\begin{enumerate}[\indent$(1)$]
\item  Les idéaux $\mathcal{H}$, $\mathcal{H}_1$ et
$\mathcal{H}_2$ représente la même classe dans $\mathbf{k}$.
\item $G=Gal(\mathbf{k}_2^{(2)}/\mathbf{k})=\langle a, b\rangle$ est un groupe
métacyclique non-modulaire  où $a^{2^n}=b^4=1$ et
$b^{-1}ab=a^{-1+k2^{n-1}}$ avec $k$ un nombre impair.
\item Seules la classe de $\mathcal{H}$ et son carré
capitulent dans chacune des trois extensions quadratiques non
ramifiées de $\mathbf{k}$.
\item Les huit classes de $C_{\mathbf{k}, 2}$ capitulent dans les trois extensions abéliennes
non ramifiées de degré $4$ de $\mathbf{k}.$
\end{enumerate}
 \end{thmp}
Dans ce qui va suivre, on adoptera les notations et les
conventions suivantes: $p$ un nombre premier, si $p\equiv 1\mod
8$, rappelons que le symbole $(\frac{2}{p})_4$ (biquadratique
rationnel) est égal à 1 ou -1, suivant que
$2^{\frac{p-1}{4}}\equiv\pm 1\mod p$. Le symbole $(\frac{p}{2})_4$
est égal à $(-1)^{\frac{p-1}{8}}$. On désigne par $h(F)$ le
$2$-nombre de classes d'un corps de nombres $F$. Rappelons aussi
que $D_3$ (resp. $Q_3$) est le groupe diédral (resp. des
quaternions) d'ordre $2^3$.
\section{Résultats préliminaires}
Ce paragraphe est réservé à certains résultats utiles dans le
reste de l'article. Les deux premiers résultats suivants
concernant des cas particuliers des extensions non ramifiées.
\bt[\cite{Hi}] \label{ta02} Soient $K/k$ une extension quadratique
et $\mu$ un nombre de $k$ premier avec $2$ tel que
$K=k(\sqrt\mu)$. L'extension de $K/k$ est non ramifiée, aux
premiers finis de $k$ si et seulement si $\mu$ vérifie les
propriétés suivantes :
\begin{enumerate}[\indent$\bullet$]
\item[$\bullet$] L'idéal principal
engendré par $\mu$ est le carré d'un idéal (fractionnaire) de $k$.
\item[$\bullet$] Il existe un nombre non nul $\xi\in k$ vérifiant
$\mu\equiv \xi^2\mod 4 $ (il s'agit d'une congruence
(multiplicative) dans $k$, modulo le sous-groupe des nombres de la
forme $1+4\frac{r}{s}$ avec $r$ et $s$ entiers de k tel que $s$ soit
premier avec $2$).
\end{enumerate}
\et

\bp[\cite{RR33}] \label{ta03} Soit $K/F$ une extension quadratique
telle que le nombre de classes de $F$ est un nombre impair. Si $K$
admet une extension non ramifiée cyclique $R$ d'ordre $4$, alors
$R/F$ est normale et $Gal(R/F)\simeq D_3$.
\ep

On aura besoin aussi de certains résultats sur les symboles
quadratiques bien connus, notamment ceux de Hilbert
$\left(\frac{a,\, b}{\mathcal{P}}\right)$ et les symboles des restes
quadratiques $[\frac{a}{\mathcal{P}}]$. (pour plus de détails, voir
\cite{Gr73} et \cite{Za99}). Explicitement la valeur du symbole de
Hilbert est donné dans un cas particulier par :

\bp[\cite{Se70}]

\label{ta04} Soient $K$ un corps de nombres de degré $n$ tel que
son anneau des entiers est principal, $a=l^{v_l(a)}u$ et
$b=l^{v_l(b)}v$ deux éléments de $K$ où $(l)$ est un idéal premier
de $K$ au-dessus d'un nombre premier $p$ tel que $p-1\geq n$ et
$v_l$ est la valuation $l$-adique. Alors
$$\left(\frac{a,\, b}{(l)}\right)=\left[\frac{-1}{(l)}\right]^{v_l(a)v_l(b)}\left[\frac{u}{(l)}\right]^{v_l(b)}\left[\frac{v}{(l)}\right]^{v_l(a)}.$$
\ep

Le résultat suivant c'est la formule de Kuroda

\bp[\cite{Lm94}]

\label{ta05} Soient $K/k$ une extension normale dont le groupe de
Galois est de type $(2, 2)$, et $k_j$ ($j=1, 2, 3$) ses trois
sous-extensions quadratiques. Alors $$h(K)=2^{d-\kappa
-2-v}q(K)h(k_1)h(k_2)h(k_3)/h(k)^2,$$ où
$q(K)=(E_K:E_{k_1}E_{k_2}E_{k_3})$ l'indice des unités de $K/k$, $d$
le nombre des premiers infinis de $k$ qui sont ramifiés dans $K$,
$\kappa $ est le $\ZZ$-rang de groupe des unités $E_k$ de $k$ et
$v=1$ ou $0$ suivant que $K\subseteq k(\sqrt{E_k})$ ou non.
\ep

Dans la suite on va rappeler des résultats concernant la théorie des
groupes qui se révélerons très utiles dans la suite de ce travail.
\bd

On dit qu'un groupe fini $G$ est métacyclique s'il possède un sous
groupe cyclique normal $H$ tel que le quotient $G/H$ est cyclique.
\ed

\bt[\cite{Hu67}]

Soit $p$ un nombre premier, un $p$-groupe métacyclique d'ordre $p^N$
peut être représenter par $$G=\langle a, b : a^{p^n}=1,\,
b^{p^m}=a^i,\, b^{-1}ab=a^q\rangle$$ avec les conditions suivantes :
\begin{enumerate}[\rm\indent(i)]
\item $n+m=N$,
\item $q^{p^m} \equiv 1\mod p^n$,
\item $i(q-1)\equiv 0\mod p^n$.
\end{enumerate}
\et

Remarquons que les groupes métacycliques abéliens sont les groupes
cycliques ou les groupes abéliens de rang 2, alors on suppose dans
toute la suite qu'un groupe métacyclique est non abélien.

\brem

Soit $G$ un groupe métacyclique, alors le groupe des commutateurs
 $G'$ est cyclique.
\erem

\bpr

Comme $G$ est un groupe métacyclique, alors il existe un sous groupe
normal $H$ de $G$ tel que $G/H$ est cyclique, donc abélien, par
suite $G'\subset N$, ce qui montre que $G'$ est cyclique.
\epr

\bd Le groupe modulaire  c'est un groupe $G$ d'ordre $ 2^n$ $(n>3)$
métacyclique tel que $G/G'\simeq (2, 2^{n-2}).$ En particulier $G'$
est d'ordre deux.
\ed

\bp[\cite{BeSn93}]

\label{ta06} Soient $G$ un $2$-groupe métacyclique non-modulaire
et $G'$ le groupe des commutateurs de $G$. Si le groupe $G/G'$ est
de type $(2, 2^m)$ avec $m>1$ et $G=\langle a, b\rangle$ tel que
$a^2\equiv b^{2^m}\mod G'$. Alors $G'=\langle a^2\rangle$ et $G$
est de l'un des types suivants :
\begin{description}
\item[{\bf Type 1}] $a^{2^{\alpha}}=1$, $b^{2^m}=1, b^{-1}ab=a^{-1}$, $\alpha>1$ ;
\item[{\bf Type 2}] $a^{2^{\alpha}}=1$, $b^{2^m}=a^{2^{\alpha -1}}$,
$b^{-1}ab=a^{-1}$, $\alpha>1$ ;
\item[{\bf Type 3}] $a^{2^{\alpha}}=1$,
$b^{2^m}=1$, $b^{-1}ab=a^{-1+k2^s}$, $1<s<\alpha$,  $k$ impair ;
\item[{\bf Type 4}] $a^{2^{\alpha}}=1$, $b^{2^m}=a^{2^{\alpha -1}}$,
$b^{-1}ab=a^{-1+k2^s}$, $1<s<\alpha$,  $k$ impair.
\end{description}
\ep

Soient $k$ un corps de nombres dont le $2$-groupe de classes est de
type $(2, 4)$ et $G=Gal(k_2^{(2)}/k)$, alors $G/G'$ est de type $(2,
4)$, donc $G=\langle a, b\rangle$ tel que $a^2\equiv b^4\equiv 1\mod
G'$ (Théorème de la base de Burnside) et $C_{k, 2}=\langle\tau,
\sigma\rangle\simeq \langle aG', bG'\rangle$ où $(\tau,
k_2^{(2)}/k)=aG'$ et $(\sigma, k_2^{(2)}/k)=bG'$ avec $(.\;,
k_2^{(2)}/k)$ est le symbole d'Artin dans $k_2^{(2)}/k$. Par suite,
il existe trois sous-groupes de $G$
d'indice 2 : $H_{1, 2}$, $H_{2, 2}$ et $H_{3, 2}$ tels que\\
\centerline{
\begin{tabular}{c}
$H_{1, 2}=\langle b, G'\rangle$, $H_{2, 2}=\langle ab, G'\rangle$ et
$H_{3, 2}=\langle a, b^2, G'\rangle$
\end{tabular}.}\\
Il existe aussi trois sous-groupes de $G$ d'indice 4 : $H_{1, 4}$,
 $H_{2, 4}$ et $H_{3, 4}$ tels que\\

\centerline{
\begin{tabular}{c}
$H_{1, 4}=\langle a, G'\rangle$, $H_{2, 4}=\langle ab^2, G'\rangle$
et $H_{3, 4}=\langle b^2, G'\rangle$
\end{tabular}.}\
\\
\indent Soient $H$ un sous-groupe de $G$ d'indice 2 ou 4, $K$ un
sous-corps de $k_2^{(2)}/k$ laissé fixe par $H$ et
$j=j_{k\rightarrow K}$ l'application de $C_{k, 2}$ vers $C_{K, 2}$,
qui fait correspondre à la classe d'un idéal $I$ de $k$ la classe de
l'idéal engendré par $I$ dans $K$.  Artin a prouvé que

\bp[\cite{Mi89}]

Il existe un homomorphisme de groupe ${ \rm V}_{G\rightarrow H}$ de
$G/G'$ vers $H/H'$ appelé le transfer de $G$ vers $H$ tel que le
diagramme suivant est commutatif

$$
\begin{CD}
   C_{k, 2} @>j>> C_{K, 2}\\
   @V(.\;, k_2^{(2)}/k)VV @VV(.\;, K_2^{(2)}/K)V\\
   G/G' @> { \rm V}_{G\rightarrow
H}>>H/H'
\end{CD}
$$
où les flèches verticales sont des isomorphismes donnés par la loi
de réciprocité d'Artin et $(.\;, k_2^{(2)}/k)$ (resp. $(.\;,
K_2^{(2)}/K)$) est le symbole d'Artin dans $k_2^{(2)}/k$ (resp.
$K_2^{(2)}/K$).

\ep

Comme les $H_{r, s}$ sont des sous-groupes normaux de
$G=Gal(k_2^{(2)}/k)$, nous utilisons la Proposition suivante
(\cite{Mi89}) pour trouver les classes de $k$ qui capitulent dans
les extensions $K_{r, s}$ ($K_{r,  s}$ est le sous-corps de
$k^{(2)}_2$ laissé fixe par $H_{r, s}$).

\bp \label{ta07} Soit $H$ un sous groupe normal d'un groupe $G$.
Pour $g\in G$, on pose $f=[\langle g\rangle.H:H]$ et soit $\{x_1,
x_2, ....,x_t\}$ un ensemble de représentants de $G/\langle g\rangle
H$; alors on a
$${\rm V}_{G\rightarrow H}(gG')=\prod_{i=1}^tx_i^{-1}g^fx_i.H'.$$
\ep

On finit par un résultat concernant le $2$-groupe de classes des
corps de nombres de type $(2^n, 2^m)$ où $n$ et $m$ sont deux
entiers strictement positifs.

\bp[\cite{BeLeSn98}]

\label{ta01}
 Soient $k$ un corps de nombres dont le $2$-groupe de
classes est de type $(2^n,  2^m)$ où $n$ et $m$ sont deux entiers
strictement
 positifs, et $k_2^{(1)}$ le $2$-corps de classes de Hilbert de $k$. S'il existe une extension quadratique non ramifiée de $k$ dont le $2$-nombre de classes
est égale à $2^{n+m-1}$, alors le $2$-nombre de classes des trois
extensions quadratiques non ramifiées de $k$ est égale à $2^{n+m-1}$
et la suite des $2$-corps de classes de Hilbert s'arrête en
$k_2^{(1)}$.

\ep
\section{Capitulation dans le corps de genres de $\mathbf{k}$}
Soient $p$ un nombre premier tel que $p\equiv 1\mod 8$,
$\mathbf{k}^*$ le corps des genres de $\mathbf{k}=\QQ(\sqrt{2p},
i)$, $h(m)$ le $2$-nombre de classes de $\QQ(\sqrt m)$ et $h(F)$ le
$2$-nombre de classes d'un corps de nombres $F$. Si
$F=\QQ(\sqrt{d_1}, \sqrt{d_2})$ est un corps biquadratique, $Q_F$
désigne l'indice du groupe engendré par les groupes des unités de
$\QQ(\sqrt{d_1})$, $\QQ(\sqrt{d_2})$ et $\QQ(\sqrt{d_1d_2})$ dans le
groupe des unités de $F$. Si $d_1=d\neq 2,\,3$ et $d_2=i$, alors
$Q_{F}$ est l'indice de Hasse de $K$ (voir page \pageref{ta17}). On
sait d'après \cite{HS82} que le nombre de classes qui capitulent
dans une extension cyclique non ramifiée $M/N$ est égal à
$[M:N][E_N:\mathcal{N}_{M/N}(E_M)]$, où $E_N$ (resp. $E_M$) est le
groupe des unités de $N$ (resp. $M$), alors pour calculer le nombre
de classes qui capitulent dans $\mathbf{k}^*/\mathbf{k}$ il faut
chercher un système fondamental d'unités (SFU) de $\mathbf{k}^*$.
Comme $\mathbf{k}^*=\QQ(\sqrt p, \sqrt 2, i)$, on va chercher un SFU
de $\mathbf{F}=\QQ(\sqrt p, \sqrt 2)$ afin de trouver un SFU de
$\mathbf{k}^*$ ( pour plus de détail sur cette méthode voir
\cite{Az99}). \bl \label{Q=2} Soient $p$, $p'$, $q_1$, $q_2$ et $q$
des nombres premiers différents tels que $p\equiv p'\equiv
-q_1\equiv -q_2 \equiv -q\mod 4$, $\pi\in\{2$, $p'$, $q$, $2q$,
$q_1q_2\}$ et $\mathbf{F}=\QQ(\sqrt p,\sqrt \pi )$. Alors l'indice
d'unités $Q_{\mathbf{F}}$ est égal à $2$. \el

\bpr D'après  ~\cite{Wa66}, on a
$h(\mathbf{F})=\frac{Q_{\mathbf{F}}h(p)h(\pi)h(p\pi)}{4}$. Or dans
tous les cas de $\pi$, on a $h(\pi)=1$ et on a aussi $h(p)=1$, alors
$h(\mathbf{F})=\frac{Q_{\mathbf{F}}h(p\pi)}{4}.$ On trouve dans
~\cite{AzMo01} que $h(\mathbf{F})=\frac{h(p\pi)}{2}$, ce qui prouve
que $Q_{\mathbf{F}}=2$. \epr

\bt \label{ta09} Soient $p$ un nombre premier tel que $p\equiv 1\mod
8$, $\mathbf{F}=\QQ(\sqrt p, \sqrt 2)$,
 $\mathbf{k}^* =\mathbf{F}(\sqrt{ -1})$ et
$\epsilon_1$ (resp. $\epsilon_2$, $\epsilon_3$) l'unité fondamentale
de $\QQ(\sqrt p)$ (resp. $\QQ(\sqrt 2)$, $\QQ(\sqrt{2p})$). Alors
\begin{enumerate}[\rm\indent(i)]
\item $\{\sqrt{\epsilon_1\epsilon_2\epsilon_3},
\epsilon_2, \epsilon_3\}$ est un SFU de $\mathbf{F}$  si et
seulement si $\epsilon_3$ est de norme $-1$.
\item $\{\epsilon_1, \epsilon_2, \sqrt{\epsilon_3}\}$ est un SFU de
$\mathbf{F}$ si et seulement si $\epsilon_3$ est de norme $1$.
\end{enumerate}
Dans les deux cas un SFU de $\mathbf{F}$ est un SFU de
$\mathbf{k}^*$.
\et

\bpr Comme l'indice des unités de $\mathbf{F}$ est égal à 2 et les
deux unités $\epsilon_1$ et $\epsilon_2$ sont de norme -1, alors
$\{\sqrt{\epsilon_1\epsilon_2\epsilon_3}$, $\epsilon_2$,
$\epsilon_3\}$ ou $\{\epsilon_1, \epsilon_2, \sqrt{\epsilon_3}\}$
est un SFU de $\mathbf{F}$ suivant que $\epsilon_3$ est de norme -1
ou 1 (\cite{Kub56}). Si $\epsilon_3$ est de norme -1, alors d'après
\cite{Az99} $\{\sqrt{\epsilon_1\epsilon_2\epsilon_3}$, $\epsilon_2$,
$\epsilon_3\}$ est un SFU de $\mathbf{k}^*$ si et seulement si il
n'existe pas d'entiers $\alpha, \beta, \gamma\in\{0, 1\}$ et qui ne
sont pas tous nuls tels que $(2+\sqrt
2)\sqrt{\epsilon_1\epsilon_2\epsilon_3}^\alpha\epsilon_2^\beta\epsilon_3^\gamma$
est un carré dans $\mathbf{F}$. Supposons que $(2+\sqrt
2)\sqrt{\epsilon_1\epsilon_2\epsilon_3}^\alpha\epsilon_2^\beta\epsilon_3^\gamma=X^2$
avec $X\in\mathbf{F}$ et les conditions précédentes. Soit $\varrho$
le $\QQ$-automorphisme défini  par $\sqrt 2\longmapsto -\sqrt 2$ et
$\sqrt p\longmapsto \sqrt p$, alors
$(X\varrho(X))^2=2\epsilon_1^\alpha(-1)^\beta(-1)^\gamma=2\epsilon_1^\alpha(-1)^{\beta+\gamma}=\pm
2\epsilon_1^\alpha$, ce qui implique que  2 est un carré dans
$\QQ(\sqrt p)$ ou bien $2\epsilon_1$ est un carré dans $\QQ(\sqrt
p)$ et ce n'est pas le cas. Si $\epsilon_3$ est de norme 1 on
reprend la même démonstration, et on trouve des contradictions. \epr
\brem

Soient $p$ un nombre premier impair, $Q_{\mathbf{k}}$ l'indice
d'unités de $\mathbf{k}$ et $\epsilon$ l'unité fondamentale de
$\QQ(\sqrt{2p})$. Alors $\epsilon$ est un carré dans $\QQ(\sqrt p,
\sqrt 2)$ si et seulement si $Q_{\mathbf{k}}=2$ si et seulement si
$\epsilon$ est de norme $1$. \erem \bpr On trouve dans \cite{Az99}
que $Q_{\mathbf{k}}=2$ si et seulement si $\epsilon$ est de norme
$1$ et d'après \cite{Az00}, on a $Q_{\mathbf{k}}=2$ si et seulement
si $2\epsilon$ est un carré dans $\mathbf{k}$. Alors pour obtenir la
remarque il suffit d'observer que $\epsilon$ est un carré dans
$\QQ(\sqrt p, \sqrt 2)$ si et seulement si $2\epsilon$ est un carré
dans $\mathbf{k}$.

\epr
\bt
\label{ta10}
Soient $\mathbf{k}=\QQ(\sqrt{2p}, i)$ avec $p$ un nombre premier tel
que $p\equiv 1\mod 8$, $\mathbf{k}^*=\QQ(\sqrt p,  \sqrt 2, i)$ le
corps de genres de $\mathbf{k}$, $C_{\mathbf{k}, 2}$ le $2$-groupe
de classes de $\mathbf{k}$ au sens large et $\epsilon$ l'unité
fondamentale de $\QQ(\sqrt{2p})$. Alors $C_{\mathbf{k}, 2}$$\simeq$
$(2^n, 2^m)$ ($n>0$ et $m>1$) et deux ou quatre classes de
$C_{\mathbf{k},  2}$ capitulent dans $\mathbf{k}^*$, suivant que
$\epsilon$ est de norme $-1$ ou $1$. \et

\bpr

Comme $p\equiv 1\mod 8$, alors d'après \cite{McPaRa95}, le
2-rang($C_{\mathbf{k}, 2})=2$, donc $C_{\mathbf{k}, 2}$ est de type
$(2^n, 2^m)$. On peut conclure facilement  que $n\geq 1$ et $m\geq
2$ (en utilisant la formule de Wada \cite{Wa66} et les résultats de
Kaplan \cite{Ka73}). Soit $\epsilon_1$ (resp. $\epsilon_2$,
$\epsilon_3$) l'unité fondamentale de $\QQ(\sqrt p)$ (resp.
$\QQ(\sqrt 2)$, $\QQ(\sqrt{2p})$) et $\mathbf{F}=\QQ(\sqrt p, \sqrt
2 )$. D'après le Théorème \ref{ta09} et \cite{Az99}, on a les
propriétés suivantes :
\begin{enumerate}[\indent$\bullet$]
\item Si $\epsilon$ est de norme -1,
alors $E_{\mathbf{k}^*}=\langle\zeta_8,
\sqrt{\epsilon_1\epsilon_2\epsilon_3}, \epsilon_2,
\epsilon_3\rangle$, ainsi
$\mathcal{N}_{\mathbf{k}^*/\mathbf{k}}(E_{\mathbf{k}^*})=\langle i,
\epsilon_3\rangle=E_{\mathbf{k}}$.
\item Si $\epsilon$ est
de norme 1, alors $E_{\mathbf{k}^*}=\langle\zeta_8, \epsilon_1,
\epsilon_2, \sqrt{\epsilon_3}\rangle$, ainsi
$\mathcal{N}_{\mathbf{k}^*/\mathbf{k}}(E_{\mathbf{k}^*})=\langle i,
\epsilon_3\rangle$, or $E_{\mathbf{k}}=\langle i,
\sqrt{i\epsilon_3}\rangle$.
\end{enumerate}

\indent Puisque le nombre de classes qui capitulent dans
$\mathbf{k}^*/\mathbf{k}$ est égal à
$2.[E_{\mathbf{k}}:\mathcal{N}_{\mathbf{k}^*/\mathbf{k}}(E_{\mathbf{k}^*})]$,
alors on a le résultat du Théorème.

\epr
\section{Commutativité de $G=Gal(\mathbf{k}^{(2)}_2/\mathbf{k})$}

Le Théorème suivant donne deux conditions nécessaires et suffisantes
pour que la tour des $2$-corps de classes de Hilbert de $\mathbf{k}$
ne s'arrête pas en premier terme, ces conditions conditions
caractérisent la Commutativité de
$G=Gal(\mathbf{k}^{(2)}_2/\mathbf{k})$ le groupe de Galois de
$\mathbf{k}^{(2)}_2/\mathbf{k}$.

\bt Soient $\mathbf{k}=\QQ(\sqrt{2p}, i)$ avec $p$ un nombre premier
tel que $p\equiv 1\mod 8$, $\mathbf{k}_2^{(1)}$ le $2$-corps de
classes de Hilbert de $\mathbf{k}$ et $\mathbf{k}_2^{(2)}$ le
$2$-corps de classes de Hilbert de $\mathbf{k}_2^{(1)}$. Alors on a
$$\begin{tabular}{c}
 $\mathbf{k}_2^{(1)}\neq\mathbf{k}_2^{(2)}$ $\Leftrightarrow $ $p=x^2+32y^2$
$\Leftrightarrow $ $(\frac{2}{p})_4=(\frac{p}{2})_4$
\end{tabular}$$
\et

\bpr

Comme $\mathbf{k}^*=\QQ(\sqrt 2, \sqrt p, i)$ est une extension de
type $(2, 2, 2)$ sur $\QQ$, alors d'après \cite{Wa66}, on a
 $$ h(\mathbf{k}^*) =
\frac{q(\mathbf{k}^*/\QQ)}{2^5}h(2)h(p)h(-1)h(-2)h(-p)h(2p)h(-2p).$$
\indent De plus $h(2)=h(p)=h(-1)=h(-2)=1,$ et
$h(\mathbf{k})=\frac{h(2p)h(-2p)}{2Q_{\mathbf{k}}}.$ Ce qui implique
que
$h(\mathbf{k^*})=\frac{q(\mathbf{k}^*/\QQ)h(-p)h(\mathbf{k})}{2^4Q_{\mathbf{k}}},$
où $q(\mathbf{k}^*/\QQ)=[E_{\mathbf{k}^*}: \langle i, \epsilon_1,
\epsilon_2, \epsilon_3\rangle],$ avec $\epsilon_1$ (resp.
$\epsilon_2$, $\epsilon_3$) l'unité fondamentale de $\QQ(\sqrt p)$
(resp. $\QQ(\sqrt 2)$, $\QQ(\sqrt{2p})$). Dans les deux cas de la
norme de $\epsilon_3$, il est facile de voir que
$q(\mathbf{k}^*/\QQ)=4$. Par suite, puisque $\mathbf{k}^*$ est une
extension  non ramifiée de $\mathbf{k}$ et le rang du 2-groupe des
classes de $\mathbf{k}$ est 2, donc
$\mathbf{k}_2^{(1)}=\mathbf{k}_2^{(2)}\Leftrightarrow
h(\mathbf{k^*})=h(\mathbf{k})/2\Leftrightarrow
h(-p)=2Q_{\mathbf{k}}$, c'est équivalent à $h(-p)=4$ et
$Q_{\mathbf{k}}=2$ ou $h(-p)=2$ et $Q_{\mathbf{k}}=1$. Or si
$h(-p)=4$ on a $Q_{\mathbf{k}}=2$ (\cite{Sc34}). D'autre part P.
Barruccand et H. Cohn ont montré dans \cite{BC69} que $h(-p)=4$ si
et seulement si $p\neq x^2+32y^2$, on déduit alors que
$\mathbf{k}_2^{(1)}\neq\mathbf{k}_2^{(2)}$ si et seulement si $p=
x^2+32y^2$. Pour completer la preuve du Théorème, on a besoin du
Lemme suivant :
\epr
\bl

Soit p un nombre premier tel que $p\equiv 1\mod 8$. Alors
$$p=x^2+32y^2\Leftrightarrow \left(\frac{2}{p}\right)_4=\left(\frac{p}{2}\right)_4$$
\el \bpr

Comme $p\equiv 1\mod 8$, alors $p=x^2+2b^2$. D'après \cite{Ka76}, on
a $(\frac{2}{p})_4(\frac{p}{2})_4=(-1)^{\frac{b}{2}}$, par suite
$p=x^2+32y^2\Leftrightarrow b=4y \Leftrightarrow
(\frac{2}{p})_4.(\frac{p}{2})_4=1$. \epr

\bc

Soient $\mathbf{k}=\QQ(\sqrt{2p},  i)$ avec p un nombre premier tel
que $p\equiv 1\mod 8$ et $C_{\mathbf{k}, 2}$ le $2$-groupe de
classes au sens large de $\mathbf{k}$. Si
$C_{\mathbf{k}, 2}$ est de type $(2, 4)$, alors \\
$$\mathbf{k}_2^{(1)}\neq\mathbf{k}_2^{(2)}\Leftrightarrow
\left(\frac{2}{p}\right)_4=\left(\frac{p}{2}\right)_4=-1.$$ \ec \bpr

Comme  $C_{\mathbf{k}, 2}$ est de type $(2, 4)$, alors d'après
\cite{Wa66} le $2$-nombre de classe $h(\mathbf{k})$ de $\mathbf{k}$
est donné par :
$$h(\mathbf{k})={\large\frac{1}{2}}Qh(2p)h(-2p)$$ o\`u $Q=Q_{\mathbf{k}}$
désigne l'indice des unités de $\mathbf{k}$, selon \cite{Ka73}
$4/h(2p)$ et $4/h(-2p)$ ou $2/h(2p)$ et $4/h(-2p)$, alors
$h(\mathbf{k})=8$ si et seulement si $h(2p)=h(-2p)=4$ et $Q=1$ ou
$h(-2p)=2h(2p)=4$ et $Q=2$. A. Scholz a montré dans \cite{Sc34} que
les derniers conditions sont équivalentes à
$(\frac{2}{p})_4=(\frac{p}{2})_4=-1$ ou bien
$(\frac{2}{p})_4=-(\frac{p}{2})_4=-1$. Alors d'après le Théorème
précédent, dans le premier cas on a
$\mathbf{k}_2^{(1)}\neq\mathbf{k}_2^{(2)}$ et dans le deuxième on a
$\mathbf{k}_2^{(1)}=\mathbf{k}_2^{(2)}$. \epr
\section{Les sous extensions de $\mathbf{k}^{(1)}_2/\mathbf{k}$}
 Dans toute cette section on suppose que $p$ est un nombre premier tel que
$p\equiv 1\mod 8$ et $(\frac{2}{p})_4=(\frac{p}{2})_4=-1$, alors
$C_{\mathbf{k}, 2}=\langle\sigma,  \tau \rangle$ où $\sigma^4=
\tau^2 $ et $\sigma \tau =\tau \sigma $, car $C_{\mathbf{k}, 2}$ est
de type $(2, 4)$. Il est clair que $C_{\mathbf{k}, 2}$ admet trois
sous groupes d'indice 2 et trois sous groupes d'indice 4. Par la
théorie du corps de classes, on sait que, chaque sous groupe $H$ de
$C_{\mathbf{k}, 2}$ correspond à une extension non ramifiée $K$ de
$\mathbf{k}_2^{(2)}$ telle que $C_{\mathbf{k}, 2}/H\simeq
Gal(K/\mathbf{k})$ et $H=\mathcal{N}_{K/\mathbf{k}}(C_{\mathbf{K},
2})$. La situation est schématisée par le diagramme suivant :
\begin{figure}[H]
$$ \xymatrix{
   & \mathbf{k}_2^{(2)} \ar@{<-}[d] & \\
   & \mathbf{k}_2^{(1)}\ar@{<-}[ld]\ar@{<-}[d]\ar@{<-}[rd]\\
  \mathbf{K_{1,  4}}\ar@{<-}[rd]&\ar@{<-}[ld] \mathbf{K_{3,  4}}\ar@{<-}[d]\ar@{<-}[rd]  & \mathbf{K_{2,  4}}\ar@{<-}[ld]\\
  \mathbf{K_{1,  2}}\ar@{<-}[rd]& \mathbf{K_{3,  2}} \ar@{<-}[d]  & \mathbf{K_{2,  2}}\ar@{<-}[ld]\\
  &\mathbf{k}
  }
$$
{\bf Diagramme 1.}
\end{figure}
Dans cette section on va essayer de construire les corps
$\mathbf{K_{1,  2}}$, $\mathbf{K_{2,  2}}$, $\mathbf{K_{3,  2}}$,
$\mathbf{K_{1,  4}}$, $\mathbf{K_{2,  4}}$, $\mathbf{K_{3,  4}}$, et
$\mathbf{k}_2^{(1)}$. Pour cela on aura besoins des deux résultats
suivants :
\brem

Si on garde les notations précédentes, alors
$$Gal(\mathbf{K_{1,  4}}/\QQ(i))\simeq Gal(\mathbf{K_{2,  4}}/\QQ(i))\simeq D_3$$
\erem

\bpr  Si on pose $F=\QQ(i)$ et $K$ le corps $\mathbf{K_{1, 4}}$ ou
$\mathbf{K_{2, 4}}$, la Proposition \ref{ta03} implique que
$Gal(\mathbf{K_{1, 4}}/\QQ(i))\simeq Gal(\mathbf{K_{2,
4}}/\QQ(i))\simeq D_3$.
\epr

\bl[\cite{Lm94}]

Soit $K/F$ une extension biquadratique telle que
$Gal(K/F)=\langle\rho, \varphi \rangle\simeq \ZZ/2\ZZ\times
\ZZ/2\ZZ$ ; si on pose $R=K(\sqrt\mu)$. Alors $R/F$  est normale si
et seulement si $\mu^{1-\rho }$ est un carré dans $K$ pour tout
$\rho\in G(K/F)$. Dans ce cas écrivons
$\mu^{1-\rho}=\alpha^2_{\rho}$, $\mu^{1-\tau}=\alpha^2_{\varphi}$ et
$\mu^{1-\rho\varphi}=\alpha^2_{\rho\varphi}$. Il est facile de voir
que $\alpha^{1+\rho}_{\rho}=\pm 1$ pour tout $\rho\in G(K/F)$; on
définit $S(\mu, K/F)=(\alpha^{1+\rho}_{\rho},
\alpha^{1+\varphi}_{\varphi},
\alpha^{1+\rho\varphi}_{\rho\varphi})$. Alors  a une permutation
près on a :

$$
\begin{array}{cc}
Gal(R/F)\simeq:\begin{cases}
\begin{array}{lcc}
(2, 2, 2)&\Leftrightarrow &
S(\mu, K/F)=(+1,  +1,  +1),\\
(2,  4)&\Leftrightarrow &
S(\mu, K/F)=(-1,  -1,  +1),\\
D_3&\Leftrightarrow
&S(\mu, K/F)=(-1,  +1,  +1),\\
Q_3&\Leftrightarrow &S(\mu, K/F)=(-1,  -1,  -1),
\end{array}
\end{cases}

\end{array}
$$
De plus $R$ est cyclique sur le corps fixé par $\langle\rho\rangle$
si et seulement si $\alpha_{\rho}^{1+\rho}=-1$, et est de type $(2,
2)$ dans le cas contraire.

\el

\bt \label{ta11} Soient $\mathbf{k}=\QQ(\sqrt{2p}, i)$ avec $p$ un
nombre premier tel que $p\equiv 1\mod 8$,
$(\frac{2}{p})_4=(\frac{p}{2})_4=-1$, $\mathbf{k}^*=\QQ(\sqrt{p},
\sqrt 2, i)$ le corps de genres de $\mathbf{k}$ et $\mathbf{K_{1,
2}}$, $\mathbf{K_{2,  2}}$, $\mathbf{K_{3,  2}}$, $\mathbf{K_{1,
4}}$, $\mathbf{K_{2,  4}}$, $\mathbf{K_{3,  4}}$, les
sous-extensions du diagramme $1$. Si on pose
$p=e^2+16f^2=x^2+32y^2=c^2-32d^2$, $\pi_1=e+4fi$, $\pi_2=e-4fi$,
$\pi_3=x+4y\sqrt{-2}$ et $\pi_4=c+4d\sqrt 2$ $(c$ et $d>0)$. Alors
\begin{enumerate}[\indent$\bullet$]
\item $\mathbf{K_{1, 2}}=\mathbf{k}(\sqrt{\pi_1})$, $\mathbf{K_{3,  2}}=\mathbf{k}^*$, $\mathbf{K_{2, 2}}=\mathbf{k}(\sqrt{\pi_2})$,
\item $\mathbf{K_{1, 4}}=\mathbf{k}^*(\sqrt{\pi_3})$, $\mathbf{K_{3, 4}}=\mathbf{k}^*(\sqrt{\pi_1})$, $\mathbf{K_{2, 4}}=\mathbf{k}^*(\sqrt{\pi_4})$,
\item $\mathbf{k}_2^{(1)}=\mathbf{k}^*(\sqrt{\pi_3}, \sqrt{\pi_4)}$.
\end{enumerate}
\et
\bpr

Comme les premiers $\pi_1$ et $\pi_2$(resp. $\pi_3$, $\pi_4$ )
sont ramifiées dans $\mathbf{k}/\QQ(i)$ (resp. dans
$\mathbf{k}^*/\QQ(\sqrt 2, i)$, $\mathbf{k}/\QQ(\sqrt 2)$) et
l'extension $\mathbf{k}^*/\mathbf{k}$ est non ramifiée. Les idéaux
engendrés par $\pi_1$ et $\pi_2$ (resp. $\pi_3$ et $\pi_4$) sont
des carrés d'idéaux de $\mathbf{k}$ (resp. $\mathbf{k}^*$).
Observons que $e$, $x$ et $c$ sont des nombres impairs, donc
$e\equiv x\equiv c\equiv \pm 1\equiv i^2\mod 4$, alors les
équations $\pi_i\equiv\xi^2 $ sont résolubles. Le Théorème
\ref{ta02} implique que les extensions $\mathbf{k}(\sqrt{\pi_1})$,
$\mathbf{k}(\sqrt{\pi_2})$, $\mathbf{k}^*(\sqrt{\pi_3})$ et
$\mathbf{k}^*(\sqrt{\pi_4})$ sont des extensions différentes  non
ramifiées de $\mathbf{k}$. Supposons que
$\mathbf{k}(\sqrt\pi_1)=\mathbf{k}^*(\sqrt\pi_2)$, alors il existe
un élément $t$ tel que $\pi_1=t^2\pi_2$, ce qui implique que
$p=t^2\pi_2^2$, et ce n'est pas le cas, car $\sqrt p\notin
\mathbf{k}$. Comme l'extension $\mathbf{k}^*/\QQ$ est normale et
$\mathbf{k}(\pi_i)/\QQ$ ($i=1, 2$) n'est pas normale, donc
$\mathbf{k}(\pi_i)\neq\mathbf{k}^*$. De la même fa\c con on montre
que $\mathbf{k}^*(\sqrt\pi_i)\neq\mathbf{k}^*(\sqrt\pi_1)$ ($i=3,
4$). Puisque $\pi_4>0$, alors le corps réel maximal de
$\mathbf{k}^*(\sqrt\pi_4)$ est $\QQ(\sqrt 2, \sqrt{\pi_4}, \sqrt
p)$ et pour $\mathbf{k}^*(\sqrt\pi_3)$ c'est  $\QQ(\sqrt 2, \sqrt
p)$, ainsi $\mathbf{k}^*(\sqrt\pi_3)\neq\mathbf{k}^*(\sqrt\pi_4)$.
Or il est facile de vérifier que, $S(\pi_3,
\mathbf{k}^*)/\QQ(i))=S(\pi_4, \mathbf{k}^*)/\QQ(i))=(-1,  +1,
+1)$, le Lemme précédent donne
$Gal(\mathbf{k}^*(\sqrt\pi_3))/\QQ(i))\simeq
Gal(\mathbf{k}^*(\sqrt\pi_4))/\QQ(i))\simeq D_3$. Ce qui achève la
preuve.

\epr

\brem On garde les notations précédentes. Alors Les deux corps
$\mathbf{K_{1,  2}}$ et $\mathbf{K_{2,  2}}$ sont conjugués, en
particulier  $h(\mathbf{K_{1,  2}})=h(\mathbf{K_{2,  2}})$. \erem
\section{le $2$-nombre de classes de $\QQ(\sqrt 2, \sqrt\pi, \sqrt p,  i )$}
Soit $p$ un nombre premier tel que $p\equiv 1\mod 8$, alors il
existe des entiers $c$ et $d$ tels que $p=c^2-32d^2$. Soient
$\pi=c+4d\sqrt 2$, $L=\QQ(\sqrt 2, \sqrt\pi, \sqrt p,  i )$ et $h$
le $2$- nombre de classes de $\QQ(\sqrt 2,
\sqrt{-\pi})=\QQ(\sqrt{-\pi})$. Dans cette section, on va calculer
$h(L)$ le $2$-nombre de classes de $L$ et $h$. Soit $E/F$ une
extension de corps de nombres tel que les anneaux des entiers de $E$
et $F$ sont principaux. Notons $[\frac{\;\,}{}]$ (resp.
$(\frac{\;\,}{})$) le symbole de reste quadratique de $E$ (resp.
$F$), $\mathcal{N}_{E/F}$ l'application norme de $E/F$ et $l$ (resp.
$\mathfrak{p}$) un nombre premier de $E$ (resp. $F)$ dont la norme
absolue est impaire, $v_l$ la valuation $l$-adique et
$\mathcal{N}_{E/F}((l))=(\mathfrak{p})^f$. Nous avons alors : \bp
Pour tout tout élément $a$ de $F$ tel que $v_l(a)=0$, on a
$$\left[\frac{a}{(l)}\right]=\left(\frac{a}{(\mathfrak{p})}\right)^f.$$
\ep \bpr Rappelons que le symbole de Hilbert sur $E$ a la Propriété
suivante :
$$\left(\frac{x,\, y}{{\beta}}\right)=\left(\frac{x,\, \mathcal{N}_{E/F}(y)}{P}\right)$$ pour $x\in F$,
$y\in E$ et $P$ un idéal premier de $F$ au-dessus de l'idéal premier
$\beta$ de $E$. Comme $v_l(a)=0$ et d'après la Proposition
\ref{ta04}, on a
$$\left[\frac{a}{(l)}\right]=\left(\frac{a,\, l}{(l)}\right)=\left(\frac{a,\, \mathcal{N}_{E/F}((l))}{(\mathfrak{p})}\right)=\left(\frac{a,\, \mathfrak{p}^f}{(\mathfrak{p})}\right)=\left(\frac{a,\, \mathfrak{p}}{(\mathfrak{p})}\right)^f=\left(\frac{a}{(\mathfrak{p})}\right)^f.$$
\epr

Avant de démontrer le Lemme suivant, rappelons que $\QQ(\sqrt{2})$
admet deux premiers infinis, $\mathcal{P}_{\infty}$ et
$\mathcal{P'}_{\infty}$. Si $u$ un élément de $\QQ(\sqrt{2})$ nous
notons $u'$ son conjugué et $s(u)=uu'|uu'|^{-1}$.

\bl On garde les notations précédentes. Alors
\begin{enumerate}[\rm\indent (i)]
\item $\left(\dfrac{-\pi,\, \epsilon_0}{(l)}\right)=\left(\dfrac{-\pi,\, \sqrt 2}{(l)}\right)=1$ pour tout nombre premier $l$ de $\QQ(\sqrt 2)$ différent de $\pi$ et de $\sqrt 2$.
\item $\left(\dfrac{-\pi,\,u}{\mathcal{P}_{\infty}}\right)=s(u)\left(\dfrac{-\pi,\,
u}{\mathcal{P'}_{\infty}}\right)$.
\item $\left(\dfrac{-\pi,\, \epsilon_0}{(\pi)}\right)=-\left(\dfrac{-\pi,\, \epsilon_0}{(\sqrt 2)}\right)=\left(\dfrac{2}{p}\right)_4\left(\dfrac{p}{2}\right)_4$.
\item $\left(\dfrac{-\pi,\, 2+\sqrt 2}{(\pi)}\right)=\left(\dfrac{-\pi,\, 2+\sqrt 2}{(\sqrt
2)}\right)=\left(\dfrac{p}{2}\right)_4$.
\end{enumerate}
\el \bpr{\rm (i)} évident, car
$v_{l}(-\pi)=v_{l}(\epsilon_0)=v_{l}(\sqrt 2)=0$.\\

{\rm (ii)} Soit $i_1$: $\sqrt 2\longmapsto \sqrt 2$ (resp. $i_2$:
$\sqrt 2\longmapsto -\sqrt 2$) le $\QQ$-plongement de $\QQ(\sqrt 2)$
dans le $\QQ(\sqrt 2)_{\mathcal{P}_{\infty}}=\RR$ (resp. $\QQ(\sqrt
2)_{\mathcal{P'}_{\infty}}=\RR$) le complété de $\QQ(\sqrt 2)$ pour
la valeur absolue associée à $\mathcal{P}_{\infty}$ (resp.
$\mathcal{P'}_{\infty}$), alors d'après \cite{Gr03}, on a : $$
  \begin{array}{l}
    \left(\dfrac{-\pi,\,u}{\mathcal{P}_{\infty}}\right)= \left\{
                                                           \begin{array}{ll}
                                                             i_1^{-1}(\left(-\pi,\,u\right)_{\mathcal{P}_{\infty}})=i_1^{-1}(1)=1, & \hbox{si $u>0$;} \\
                                                             i_1^{-1}(\left(-\pi,\,u\right)_{\mathcal{P}_{\infty}})=i_1^{-1}(-1)=-1, & \hbox{si $u<0$.}
                                                           \end{array}
                                                         \right.
\\
    \left(\dfrac{-\pi,\,u}{\mathcal{P'}_{\infty}}\right)=\left\{
                                                          \begin{array}{ll}
                                                             i_2^{-1}(\left(-\pi',\,u'\right)_{\mathcal{P'}_{\infty}})=i_2^{-1}(1)=1, & \hbox{si $u'>0$;} \\
                                                             i_2^{-1}(\left(-\pi',\,u'\right)_{\mathcal{P'}_{\infty}})=i_2^{-1}(-1)=-1, & \hbox{si $u'<0$.}
                                                           \end{array}
                                                        \right.
  \end{array}
$$
Car $\left(v,\,u\right)_{\RR}=-1$ si et seulement si les deux
nombres $u$ et $v$ sont négatifs, où
$\left(-\pi,\,u\right)_{\mathcal{P}_{\infty}}$ (resp.
$\left(-\pi,\,u\right)_{\mathcal{P'}_{\infty}}$) est le symbole
local de Hilbert définit sur $\QQ(\sqrt
2)_{\mathcal{P}_{\infty}}\times\QQ(\sqrt 2)_{\mathcal{P}_{\infty}}$
(resp. $\QQ(\sqrt 2)_{\mathcal{P'}_{\infty}}\times\QQ(\sqrt
2)_{\mathcal{P'}_{\infty}}$). D'où le résultat.
\\

{\rm (iii)} Comme $v_{l}(-\pi)=1$ et $v_{l}(\epsilon)=0$ et d'après
la Proposition \ref{ta04}, on a $$\left(\frac{-\pi,\,
\epsilon_0}{(l)}\right)=\left[\frac{\epsilon_0}{(\pi)}\right],$$ où
$[\frac{\;\,}{}]$ est le symbole de reste quadratique sur $\QQ(\sqrt
2)$. Or $\epsilon_0=1+\sqrt 2$ et $\pi =c+4d\sqrt 2$, alors
$$\left[\frac{\epsilon_0}{(\pi)}\right]=\left[\frac{4d}{(\pi)}\right]\left[\frac{4d+4\sqrt 2d}{(\pi)}\right]=\left[\frac{4d}{(\pi)}\right]\left[\frac{4d-c+\pi}{(\pi)}\right]=\left[\frac{d}{(\pi)}\right]\left[\frac{4d-c}{(\pi)}\right].$$
D'après \cite{Ka76}, Théorème 2, page 12 et la Proposition
précédente on a
$$\left[\frac{\epsilon_0}{(\pi)}\right]=\left[\frac{d}{(p)}\right]\left[\frac{4d-c}{(p)}\right]=\left(\frac{2}{p}\right)_4\left(\frac{p}{2}\right)_4.$$
Notons $S$ l'ensemble des nombres premiers $l$ de $\QQ(\sqrt 2)$
différents de $\pi$, de $\sqrt 2$, de $\mathcal{P}_{\infty}$ et de
$\mathcal{P'}_{\infty}$ les deux premiers infinis de $\QQ(\sqrt 2)$,
alors d'après, la formule du produit pour le symbole de Hilbert, on
a $$\prod_{l\in S}\left(\frac{-\pi,\,
\epsilon_0}{(l)}\right)\left(\frac{-\pi,\,
\epsilon_0}{\mathcal{P}_{\infty}}\right)\left(\frac{-\pi,\,
\epsilon_0}{\mathcal{P'}_{\infty}}\right)\left(\frac{-\pi,\,
\epsilon_0}{(\pi)}\right)\left(\frac{-\pi,\, \epsilon_0}{(\sqrt
2)}\right)=1.$$ Puisque $s(\epsilon_0)=-1$, alors
$\left(\dfrac{-\pi,
\epsilon_0}{\mathcal{P}_{\infty}}\right)=-\left(\dfrac{-\pi,
\epsilon_0}{\mathcal{P'}_{\infty}}\right)$. Il en résulte de {\rm
(i)} que $$\left(\dfrac{-\pi,\, \epsilon_0}{(\pi)}\right)=-\left(\dfrac{-\pi,\, \epsilon_0}{(\sqrt 2)}\right)=\left(\dfrac{2}{p}\right)_4\left(\dfrac{p}{2}\right)_4$$.\\

{\rm (iv)} De la même façon que dans {\rm (iii)}, on trouve le
résultat annoncé. \epr \bt

\label{ta12} Si on garde les notations et hypothèses précédentes.
Alors le $2$-groupe de classes de $\QQ(\sqrt{-\pi})$ est cyclique.
De plus on a
$$
\begin{array}{ccc}
h=\begin{cases}
\begin{array}{rll}
1&\text{\rm si }&(\frac{2}{p})_4=-(\frac{p}{2})_4,\\
2&\text{\rm si }&(\frac{2}{p})_4=(\frac{p}{2})_4=-1,\\
\geq 4&\text{\rm si}&(\frac{2}{p})_4=(\frac{p}{2})_4=1.
\end{array}
\end{cases}
\end{array}
$$\et
\bpr On a $\QQ(\sqrt{-\pi})=\QQ(\sqrt 2)(\sqrt{-\pi})$. Le nombre
de classes de $\QQ(\sqrt 2)$ est égal à $1$, alors d'après
\cite{Gr73}, le $2$-rang du $2$-groupe de classes de
$\QQ(\sqrt{-\pi})$ est $r_2=t-1-e$, où $t$ est le nombre des
premiers de $\QQ(\sqrt{2})$ qui sont ramifiés dans
$\QQ(\sqrt{-\pi})$ et $e$ est l'entiers naturel tel que $2^e$ est
le l'indice du groupe engendré par les unités de $\QQ(\sqrt{2})$
qui sont des normes dans $\QQ(\sqrt{-\pi})$ dans le groupe des
unité de $\QQ(\sqrt{2})$. Observons que $c$ est un nombre impair
tel que $-(c+4d\sqrt 2)=-\pi\equiv -(\frac{-1}{c})\equiv
-(\frac{2}{p})_4(\frac{p}{2})_4\mod 4$ (voir \cite{Ka76}, Théorème
2, page 12) ; donc d'après la Proposition 1.2 de \cite{Gr73}, page
11, $\QQ(\sqrt{-\pi})/\QQ(\sqrt{2})$ est non ramifié en $\sqrt 2$
si et seulement si $(\frac{2}{p})_4=-(\frac{p}{2})_4$. De plus
$\QQ(\sqrt{2})$ admet deux premiers infinis,
$\mathcal{P}_{\infty}$ et $\mathcal{P'}_{\infty}$, se ramifient
dans $\QQ(\sqrt{-\pi})/\QQ(\sqrt{2})$. Comme $-\pi$ et $-1$ sont
deux nombres négatifs et $s(\epsilon_0)=-1$, alors le (ii) du
Lemme précédent donne que $\left(\dfrac{-\pi,
-1}{\mathcal{P}_{\infty}}\right)=-1$ et $\left(\dfrac{-\pi,
\epsilon_0}{\mathcal{P}_{\infty}}\right)=-\left(\dfrac{-\pi,
\epsilon_0}{\mathcal{P'}_{\infty}}\right)$, alors nous avons
toujours que $e=2$, car $-1$ et $\epsilon_0$ ne sont pas des
normes dans l'extension $\QQ(\sqrt{-\pi})/\QQ(\sqrt{2})$ (un
élément $u$ de $\QQ(\sqrt{2})$ est norme dans cette extension si
et seulement si la valeur du symbole de Hilbert
$\left(\dfrac{-\pi,\, u}{(l)}\right)=1$ pour tout idéal premier
$(l)$ de $\QQ(\sqrt 2)$). Si $(\frac{2}{p})_4=-(\frac{p}{2})_4$,
alors $e=2$ et $t=3$, par suite $h=1$. Supposons dans toute la
suite que $(\frac{2}{p})_4=(\frac{p}{2})_4$, on vérifie facilement
que $t=4$ et nous avons $e=2$, ce qui prouve que $r_2=1$ et le
$2$-groupe de classes de $\QQ(\sqrt{-\pi})$ est cyclique. Soient
$r_4$ le $4$-rang du $2$-groupe de classes de $\QQ(\sqrt{-\pi})$
et $\mathfrak{a}$ l'idéal premier de $\QQ(\sqrt{-\pi})$ tel que
$(\sqrt 2)=\mathfrak{a}^2$, il est facile de voir que la classe de
$\mathfrak{a}$ dans $\QQ(\sqrt{-\pi})$ est une classe ambigue non
trivial et l'idéal de $\QQ(\sqrt 2)$ engendré par $\sqrt 2$ est
engendré aussi par $(2+\sqrt 2)$ car $-\sqrt 2=(2+\sqrt 2)(1-\sqrt
2)=(2+\sqrt 2)\epsilon_0^{-1}$. D'après la théorie du genres,
$r_4\geq 1$ si  $(2+\sqrt 2)$ est norme dans l'extension
$\QQ(\sqrt{-\pi})/\QQ(\sqrt{2})$, ainsi le (iv) du Lemme précédent
entraîne que $r_4\geq 1$ si et seulement si
$(\frac{2}{p})_4=(\frac{p}{2})_4=1$. D'où le résultat énoncé.\epr

\label{ta17} Rappelons qu'une extension CM est une extension
quadratique totalement imaginaire d'un corps de nombres totalement
réel. Soit $K/k$ une extension CM, alors l'indice des unités de
Hasse est défini par $Q_K=[E_K:W_KE_k]$ où $W_K $ est le groupe des
racines de l'unité contenues dans $K$, $E_K$ (resp. $E_k$) le groupe
des unités de $K$ (resp. k) et on note par $\omega_K$ le cardinal de
$W_K$. Il est à noter que $Q_K=1$ ou $2$ (voir \cite{Ha85}). Dans le
Lemme suivant on va calculer l'indice des unités de Hasse de
 $M=\QQ(\sqrt 2, \sqrt{-\pi}, \sqrt p )$, $L=M(i)$ et  $\mathbf{k}^*$
le corps de genres de $\mathbf{k}=\QQ(\sqrt{2p},  i).$
\bl Si on garde les notations et hypothèses précédentes. Alors
$$Q_M=Q_{\mathbf{k}^*}=Q_L=1,$$
\el
\bpr

Soient $p$ un nombre premier tel que $p\equiv 1\mod 8$,
$\mathbf{F}=\QQ(\sqrt p, \sqrt 2)$ et $\epsilon_1$ (resp.
$\epsilon_2$, $\epsilon_3$) l'unité fondamentale de $\QQ(\sqrt p)$
(resp. $\QQ(\sqrt 2)$, $\QQ(\sqrt{2p})$). Alors
$\{\sqrt{\epsilon_1\epsilon_2\epsilon_3}$, $\epsilon_2$,
$\epsilon_3\}$ ou $\{\epsilon_1, \epsilon_2, \sqrt{\epsilon_3}\}$
est un SFU de $\mathbf{F}$ suivant que $\epsilon_3$ est de norme -1
ou 1. Si $\epsilon_3$ est de norme -1, alors d'après \cite{Az99},
$\{\sqrt{\epsilon_1\epsilon_2\epsilon_3}$, $\epsilon_2$,
$\epsilon_3\}$ est un SFU de $M$ si et seulement s'il n'existe pas
d'entiers $\alpha, \beta, \gamma\in\{0, 1\}$  qui ne sont pas tous
nuls et tels que
$\pi\sqrt{\epsilon_1\epsilon_2\epsilon_3}^\alpha\epsilon_2^\beta\epsilon_3^\gamma$
est un carré dans $\mathbf{F}$. Supposons que
$\pi\sqrt{\epsilon_1\epsilon_2\epsilon_3}^\alpha\epsilon_2^\beta\epsilon_3^\gamma$
est un carré dans $\mathbf{F}$. Comme $\pi$ (resp. $\epsilon_1,
\epsilon_2$) est de norme $p$ (resp. -1) dans
$\mathbf{F}/\QQ(\sqrt{2p})$, alors
$p\epsilon_3^\alpha(-1)^\beta\epsilon_3^{2\gamma}$ est un carré dans
$\QQ(\sqrt{2p})$, ce qui implique que $p$ est un carré dans
$\QQ(\sqrt{ 2p})$ ou bien $p\epsilon_3$ est un carré dans
$\QQ(\sqrt{2p})$ et ce n'est pas le cas. Si $\epsilon_3$ est de
norme 1 on reprend la même démonstration, et on trouve des
contradictions. Il reste de prouver que $Q_{\mathbf{k}^*}=Q_L=1$, le
Théorème \ref{ta09} implique que $Q_{\mathbf{k}^*}=1$. On en déduit,
avec le Lemme 25 de \cite{Ok01}  que $Q_L=1$.\epr

\bt

\label{ta13} Soient $L=\QQ(\sqrt 2, \sqrt\pi, \sqrt p,  i )$ avec
$p=c^2-32d^2$ un nombre premier tels que $p\equiv 1\mod 8$,
$\pi=c+4d\sqrt 2$ $(c$, $d>0)$ et $h$ le $2$-nombre de classes de
$\QQ(\sqrt{-\pi})$. Alors $$
\begin{array}{cc}
h(L)=\begin{cases}
\begin{array}{ll}
 {(\frac{h}{2})}^2.\frac{h(\mathbf{k}^*)}{2}&\text{\rm si $(\frac{2}{p})_4=(\frac{p}{2})_4=1$,}\\ \frac{h(\mathbf{k}^*)}{2}&\text{\rm sinon.
}
\end{array}
\end{cases}
\end{array}$$
\et
\bpr

Notons que $L/\mathbf{F}$ est une extension normale de type $(2, 2)$
qui vérifie les hypothèses de la Proposition 2 de \cite{Lm95}. Il en
résulte que $$
h(L)=\frac{Q_L}{Q_MQ_{\mathbf{k}^*}}\frac{\omega_L}{\omega_M\omega_{\mathbf{k}^*}}\frac{h(M)h(\mathbf{k}^*)h(N)}{h(\mathbf{F})^2},
$$ avec $M=\QQ(\sqrt p, \sqrt{-\pi})$, $N=\QQ(\sqrt
p, \sqrt{\pi})$, $\mathbf{k}^*=\QQ(\sqrt 2, \sqrt p, i)$ le corps de
genres de $\mathbf{k}=\QQ(\sqrt{2p}, i)$ et $\mathbf{F}=\QQ(\sqrt 2,
\sqrt{p})$. Lorsque $(\frac{2}{p})_4=-(\frac{p}{2})_4$, on a vu aux
sections 4 et 5 que la suite des 2-corps de classes de $\mathbf{k}$
de Hilbert s'arrête en $\mathbf{k}_2^{(1)}$ et $L$ est une extension
quadratique non ramifiée de $\mathbf{k}^*$, par suite on a
$$h(L)=\frac{h(\mathbf{k}^*)}{2}.$$
\indent Si $(\frac{2}{p})_4=(\frac{p}{2})_4$, alors
$(\frac{2}{p})_4(\frac{p}{2})_4=(\frac{-1}{c})=1$ (voir
\cite{Ka76}). Dans ce cas, $c\equiv 1\mod 4$, ainsi $\pi\equiv 1\mod
4$. Sous ces conditions, le Théorème \ref{ta02} implique que
$N=\mathbf{F}(\sqrt\pi)$  est une extension quadratique non
ramifiée, et donc puisque le groupe de classes de $\mathbf{F}$ est
cyclique (voir \cite{AzMo01}) on a $$ h(N) =\frac{h(\mathbf{F})}{2}.
\eqno{(6.1)} $$ \indent Maintenant, remarquons que $M/\QQ(\sqrt 2)$
est une extension normale de type $(2, 2)$, nous pouvons alors
appliquer la formule de Kuroda (Proposition \ref{ta05}) et on a
$$h(M)=\frac{qhh'h(\mathbf{F})}{2h(\QQ(\sqrt2))^2},$$ avec
$q=[E_M:EE'E_{\mathbf{F}}]$, $h$ (resp. $h'$) le $2$-nombre de
classes de $\QQ(\sqrt{-\pi})$ (resp. $\QQ(\sqrt{-\pi'})$),
$\pi'=c-4d\sqrt 2$ et $E$ (resp. $E'$) le groupe des unités de
$\QQ(\sqrt{-\pi})$ (resp. $\QQ(\sqrt{-\pi'})$). Comme
$\QQ(\sqrt{-\pi})$ et $\QQ(\sqrt{-\pi'})$)  sont deux corps de
nombres conjugués, alors $h=h'$ et le Lemme précédent implique que
$$[E_M:E_{\mathbf{F}}]=1.$$
\indent Il est clair que $\QQ(\sqrt{-\pi})$ et $\QQ(\sqrt{-\pi'})$
sont des extensions CM, tel que leurs indices des unités sont égaux
à 1 (Il suffit de remarquer que ${-\pi}\epsilon_0$ n'est jamais un
carré dans $\QQ(\sqrt 2)$ où $\epsilon_0$ est l'unité fondamentale
de $\QQ(\sqrt 2)$), par suite $E$ et $E'$ qui sont engendrés par -1
et $\epsilon_0$ et inclus dans $E_{\mathbf{F}}$ c'est-à-dire
$q=[E_M:EE'E_{\mathbf{F}}]=[E_M:E_{\mathbf{F}}]=1$. Ainsi
$$h(M)=\frac{1}{2}h^2h(\mathbf{F}). \eqno{(6.2)}$$
\indent Compte tenu du Lemme précédent on a
$$\frac{Q_L}{Q_MQ_{\mathbf{k}^*}}\frac{\omega_L}{\omega_M\omega_{\mathbf{k}^*}}=\frac{1}{1.1}\frac{8}{2.8}=\frac{1}{2}~.\eqno{(6.3)}$$
\indent Les résultats (6.1), (6.2) et (6.3) implique que
$h(L)={(\frac{h}{2})}^2.\frac{h(\mathbf{k}^*)}{2}$. Les résultats du
Théorème \ref{ta12} achèvent la preuve.
\epr

\section{Structure de  $G=Gal(\mathbf{k}^{(2)}_2/\mathbf{k})$ et de $C_{2, \mathbf{k}}$ }

\bt

\label{ta14} Soient $\mathbf{k}=\QQ(\sqrt{2p}, i)$ avec $p$ un
nombre premier tel que $p\equiv 1\mod 8$,
$\mathbf{k}^*=\QQ(\sqrt{2}, \sqrt{p}, i)$ le corps de genres de
$\mathbf{k}$ et $C_{\mathbf{k}^*, 2}$ le $2$-groupe de classes de
$\mathbf{k}^*$. Alors le rang de $C_{\mathbf{k}^*, 2}=2$ ou $3$. De
plus  le rang de $C_{\mathbf{k}^*, 2}=3$ si et seulement si $(
\frac{2}{p} )_4=( \frac{p}{2})_4=1$. \et \bpr Notons $F$ le corps
$\QQ(\sqrt 2, i)=\QQ(\zeta_8)$, $\mathrm{Am}(\mathbf{k}^*/F)$ le groupe de
classes ambigues dans $\mathbf{k}^*/F$ et $r$ le rang de
$C_{\mathbf{k}^*, 2}$. Il est bien connu que le groupe des unités de
$F$ est engendré par $\epsilon_0=1+\sqrt 2$ l'unité fondamentale de
$\QQ(\sqrt 2)$ et $\zeta_8$ la racine $8$-ième de l'unité, de plus
le nombre de classes de $F$ est égal à $1$. Alors la formule de
genres donne le nombre des
classes ambigues dans $\mathbf{k}^*/F$:\\
$$|Am(\mathbf{k}^*/F)|=\frac{2^3}{[E_F:E_F\cap
\mathcal{N}_{\mathbf{k}^*/F}(\mathbf{k}^*)]}=2^r,$$ car il existe
quatre idéaux premiers de $F$ qui se ramifient dans $\mathbf{k}^*$,
ces idéaux sont au-dessus de $p$. Comme $F$ est imaginaire,
$\mathbf{k}^*=F(\sqrt p)$ et grace à la formule de produit pour le
symbole  de  Hilbert $\left(\dfrac{\bullet,\,\bullet}{\beta}\right)
$; le théorème de Hasse entraîne qu'une unité $\epsilon $ de $F$ est
une norme si et seulement si $(p,\;\epsilon)_{\beta}=1$ pour tout
idéal premier de $F$  qui n'est pas au-dessus de $2$. En utilisant
les propriétés du symbole de  Hilbert et le même raisonnement que
dans le Lemme 4 on trouve que
$$\begin{array}{cc} \left(\dfrac{p,\,\epsilon}{\beta}\right)=\begin{cases}
\begin{array}{lll}
1&\text{\rm si }&\text{\rm$\beta$ n'est pas au-dessus de $p$}, \\
(\frac{p}{2})_4&\text{\rm si }&\text{\rm$\beta$ au-dessus de $p$
et $\epsilon=\zeta _8$,} \\
(\frac{p}{2})_4(\frac{2}{p})_4&\text{\rm si }&\text{\rm$\beta$
au-dessus de $p$ et  $\epsilon=\epsilon_0$,}
\\ (\frac{2}{p})_4&\text{\rm si }&\text{\rm$\beta$ au-dessus de $p$
et  $\epsilon=\epsilon_0\zeta _8$.}
\end{array}
\end{cases}
\end{array}
$$ En particulier

$$
\begin{array}{cc}
E_F\cap \mathcal{N}_{\mathbf{k}^*/F}(\mathbf{k}^*)=\begin{cases}
\begin{array}{lll}
\langle\zeta _8, \epsilon_0\rangle&\text{\rm si }&(\frac{2}{p})_4=(\frac{p}{2})_4=1,\\
\langle i, \epsilon_0\rangle&\text{\rm si }&(\frac{2}{p})_4=(\frac{p}{2})_4=-1,\\
\langle\zeta _8, \epsilon_0^2\rangle&\text{\rm si }&(\frac{2}{p})_4=-(\frac{p}{2})_4=-1,\\
\langle i, \epsilon_0\zeta _8\rangle&\text{\rm si
}&(\frac{2}{p})_4=-(\frac{p}{2})_4=1.\
\end{array}
\end{cases}
\end{array}
$$
Enfin $$
\begin{array}{cc}
|Am(\mathbf{k}^*/F)|=\begin{cases}
\begin{array}{lll}
2^3&\text{\rm si }&(\frac{2}{p})_4=(\frac{p}{2})_4=1,\\
2^2&\text{\rm sinon.}&
\end{array}
\end{cases}
\end{array}
$$ D'où le résultat. \epr \brem Soient p un nombre
premier tel que $p\equiv 1\mod 8$ et $\mathbf{k}^*$ le corps de
genres de $\mathbf{k}=\QQ(\sqrt{2p}, i)$. Alors on a les
décompositions suivantes :\\$(p)=(\pi\pi')$ dans $\QQ(\sqrt 2),$\\
$(\pi)=(\pi_1\pi_2)$ et $(\pi')=(\pi_3\pi_4)$ dans $\QQ(\zeta
_8)$,\\ $(\pi_i)=\mathcal{G}_i^2$ dans $\mathbf{k}^*$,\\
$(p)=\beta^2$ dans $\QQ(\sqrt{2p})$,\\
$\beta=\mathcal{P}_1\mathcal{P}_2$ dans $\QQ(\sqrt 2, \sqrt p)$,\\
$\mathcal{P}_1=\mathcal{G}_1\mathcal{G}_2$ dans $\mathbf{k}^*$,\\
$(\pi)=\mathcal{P}_1^2$ dans $\QQ(\sqrt 2, \sqrt p)$.
\erem
\bt

\label{ta15} Soient $\mathbf{k}=\QQ(\sqrt{2p}, i)$ avec $p$ un
nombre premier tel que $p\equiv 1\mod 8$,
$(\frac{2}{p})_4=(\frac{p}{2})_4=-1$, $\mathbf{k}^*=\QQ(\sqrt{2},
\sqrt{p}, i)$ le corps de genres de $\mathbf{k}$ et
$C_{\mathbf{k}^*, 2}$ le $2$-groupe de classes de $\mathbf{k}^*$.
Alors le $4$-rang de $C_{\mathbf{k}^*, 2}=1$ . \et \bpr D'après le
Théorème précédent, la condition
$(\frac{2}{p})_4=(\frac{p}{2})_4=-1$ entraîne que le rang de
$C_{\mathbf{k}^*, 2}$ est égal à $2$ ou encore il y a
exactement $4$ classes ambigues dans $\mathbf{k}^*/\QQ(\zeta _8)$.\\
\indent Il est à noter que, si $E/F$ est une extension quadratique
de corps de nombres telle que le nombre de classes de $F$ est impair
et il existe exactement $2^r$ classes ambigues qui sont des carrés,
alors le $4$-rang de $C_E$ est égal à $r$. Comme le $2$-nombre de
classes de $\mathbf{k}^*$ est égal à $2h(-p)$ et est divisible par
$16$, alors le $4$-rang de $C_{\mathbf{k}^*, 2}=1$ ou $2$. Il reste
de trouver la classe ambigue non triviale de $C_{\mathbf{k}^*, 2}$
qui n'est pas un carré.
\\
\indent Nous reprenons les notations de la remarque précédente.
Alors l'idéal $\beta$ est non principal car sinon $p\epsilon$ est un
carré dans $\QQ(\sqrt{2p})$ où $\epsilon$ est l'unité fondamentale
de $\QQ(\sqrt{2p})$, ce qui implique que $\epsilon$ est un carré
dans $\QQ(\sqrt 2, \sqrt{p})$ c'est-à-dire que $\epsilon$ est de
norme $1$ (voir remarque $2$). Mais puisque
$(\frac{2}{p})_4=(\frac{p}{2})_4=-1$, $\epsilon$ est de norme $-1$,
ainsi on obtient une contradiction. Or la relation
$\mathcal{N}_{\QQ(\sqrt 2,
\sqrt{p})/\QQ(\sqrt{2p})}(\mathcal{P}_1)=\beta$ implique que
$\mathcal{P}_1$ est non principal par suite la classe de
$\mathcal{P}_1$ est d'ordre $2$. De même  on a
$\mathcal{N}_{\mathbf{k}^*/\QQ(\sqrt 2,
\sqrt{p})}(\mathcal{G}_1)=\mathcal{P}_1$ et $\mathcal{G}_1$ est non
principal dans $\mathbf{k}^*$, alors la classe $[\mathcal{G}_1]$ est
ambigue non triviale dans $\mathbf{k}^*/\QQ(\zeta _8)$ car
$(\pi_1)=\mathcal{G}_1^2$. Rappelons que puisque
$(\frac{2}{p})_4=(\frac{p}{2})_4=-1$, alors le nombre de classes de
$\QQ(\sqrt 2, \sqrt{p})$ est égal à $2$, ce qui prouve que la classe
$[\mathcal{P}_1]=\mathcal{N}_{\mathbf{k}^*/\QQ(\sqrt
2\sqrt{p})}([\mathcal{G}_1])$ engendre le groupe de classes de
$\QQ(\sqrt 2, \sqrt{p})$, ainsi la classe $[\mathcal{G}_1]$ n'est
pas un carré dans $C_{\mathbf{k}^*, 2}$. D'où le résultat
énoncé.\epr \bl Soit $L/k$ une extension quadratique non ramifiée
telle que la suite des $2$-corps de Hilbert de $k$ s'arrête en
$k_2^{(1)}$. Alors $\mathcal{N}_{L/k}(C_{L, 2})\simeq C_{L, 2}$. \el
\bpr Montrons que l'homomorphisme suivant est injectif
$$\begin{array}{ccc}

C_{L, 2} &\longrightarrow  &C_{k, 2}\\
c &\longmapsto  &\mathcal{N}_{L/k}(c)
\end{array}$$
Comme $L$ est une extension non ramifiée de $k$, alors d'après
la théorie du corps de classes on a\\
$$[C_{k, 2}:\mathcal{N}_{L/k}(C_{L, 2})]=[L:k]=2,$$ ce qui implique que
$$|\mathcal{N}_{L/k}(C_{L, 2})|=\frac{h(k)}{2}.$$ Puisque que la suite
des $2$-corps de Hilbert de $k$ s'arrête en $k_2^{(1)}$ et que $L$
est une extension quadratique non ramifiée de $k$, alors
$$h(L)=\frac{h(k)}{2}.$$ Ce qui prouve que
$h(L)=|\mathcal{N}_{L/k}(C_{L, 2})|$, par suite $\mathcal{N}_{L/k}$
est injectif et
$$\mathcal{N}_{L/k}(C_{L, 2})\simeq C_{L, 2}.$$
\epr
\bt

\label{ta16} Soient $\mathbf{k}=\QQ(\sqrt{2p}, i)$ avec
$p=c^2-32d^2$ un nombre premier tel que $p\equiv 1\mod 8$,
$\pi=c+4d\sqrt 2>0$, $(\frac{2}{p})_4=(\frac{p}{2})_4=-1$,
$\mathbf{k}_2^{(1)}$ le $2$-corps de classes de Hilbert de
$\mathbf{k}$, $\mathbf{k}_2^{(2)}$ le $2$-corps de classes de
Hilbert de $\mathbf{k}_2^{(1)}$ et
$G=Gal(\mathbf{k}_2^{(2)}/\mathbf{k})$. Alors $G$ est métacyclique
non-modulaire et la suite des $2$-corps de classes de Hilbert de
$\mathbf{k}$ (resp. $\mathbf{k}^*$) s'arrête en $\mathbf{k}_2^{(2)}$
(resp. ${\mathbf{k}^*}_2^{(1)}=\mathbf{k}_2^{(2)}$). \et \bpr Compte
tenu des Théorèmes \ref{ta14} et \ref{ta15}, on peut conclure que le
2-groupe de classes de $\mathbf{k}^*$ est de type $(2, 2^m)$, alors
$\mathbf{k}^*$ admet trois extensions quadratiques non ramifiées. Le
diagramme 1 et le Théorème \ref{ta11}, montrent que ces trois
extensions non ramifiées sont : $\mathbf{K_{1, 4}}$, $\mathbf{K_{2,
4}}$ et $\mathbf{K_{3, 4}}$ et on a donc, d'après le Théorème
\ref{ta13}, $h(\mathbf{K_{2, 4}})=\frac{h(\mathbf{k}^*)}{2}=h(-p)$.
Dans ce cas, la Proposition \ref{ta01} implique que la suite des
$2$-corps de classes de Hilbert de $\mathbf{k}^*$ s'arrête en
${\mathbf{k}^*}_2^{(1)}$. Par ailleurs, il découle du lemme
précédent que $\mathcal{N}_{M_i/\mathbf{k}}(C_{M_i, 2})\simeq
C_{M_{i, 2}}$, où $M_i=\mathbf{K}_{i, 4}$. La théorie du corps de
classes nous donne que $\mathcal{N}_{M_i/\mathbf{k}}(C_{M_i, 2})$
est cyclique pour deux indice $i$. On en déduit, que $M_i$ et
$\mathbf{k}_2^{(1)}$ ont le même $2$-corps de classes de Hilbert
$\mathbf{k}_2^{(2)}$, donc $G'$ est d'ordre $h(-p)/2\geq 4$. Par
suite ${\mathbf{k}^*}_2^{(1)}=\mathbf{k}_2^{(2)}$ et $G$ est
non-modulaire. Soit $i$ tel que $M_i/\mathbf{k}$ est cyclique, alors
$H=Gal(\mathbf{k}_2^{(2)}/M_i)$ est un sous-groupe cyclique normal
de $G=Gal(\mathbf{k}_2^{(2)}/\mathbf{k})$ tel que $G/H\simeq
Gal(M_i/\mathbf{k})$. Alors $G$ est un groupe métacyclique
non-modulaire, ce qui entraîne que $G'$ est cyclique, par conséquent
 $\mathbf{k}_2^{(2)}=\mathbf{k}_2^{(3)}$.\epr
 \bl
Soit $L/k$ une extension quadratique ramifiée. Alors l'homomorphisme
suivant est surjectif $$\begin{array}{ccc}

C_{L, 2} &\longrightarrow  &C_{k, 2}\\
c &\longmapsto  &\mathcal{N}_{L/k}(c)
\end{array}$$
\el
\bpr

Comme $L/k$ est ramifiée, la théorie du corps de classes implique
que
$$[C_{k, 2}:\mathcal{N}_{L/k}(C_{L, 2})]<[L:k]=2.$$
Autrement dit $\mathcal{N}_{L/k}$ est surjective.
\epr
\bp

\label{ta08} Soient $d$ un entier naturel sans facteurs carrés,
$k=\QQ(\sqrt{2d}, i)$, $\epsilon$ l'unité fondamentale de $k$ et
$\mathcal{H}$ l'idéal premier au-dessus de $1+i$ dans $k$. Si
l'indice des unités de $k$ est égal à $1$, alors la classe de
$\mathcal{H}$ dans $k$ est d'ordre $2$. De plus la classe
$\mathcal{H}$ capitule dans $k(\sqrt 2)$. \ep \bpr même
démonstration qui se trouve dans \cite{Az00}. \epr \bt Soient
$\mathbf{k}=\QQ(\sqrt{2p}, i)$ avec $p$ un nombre premier tel que
$p\equiv 1\mod 8$, $(\frac{2}{p})_4=(\frac{p}{2})_4=-1$,
$C_{\mathbf{k}, 2}$ le $2$-groupe de classes de $\mathbf{k}$ et
$\mathcal{H}$ l'idéal premier au-dessus de $1+i$ dans $\mathbf{k}$.
Alors $C_{\mathbf{k}, 2}=\langle\sigma, \tau\rangle$ tel que
$\sigma^2=[\mathcal{H}] $, la classe de $\mathcal{H}$ dans
$\mathbf{k}$ et $\mathcal{N}_{\mathbf{k}/\QQ(\sqrt{2p})}(\tau)=1$.
\et \bpr

Comme $(\frac{2}{p})_4=(\frac{p}{2})_4=-1$, alors le $2$-groupe de
classes de $\QQ(\sqrt{2p})$ est cyclique d'ordre $4$, engendré par
une classe $c$ et $$C_{\mathbf{k}, 2}=\langle\sigma, \tau\rangle
\text{ tel que } \sigma^4=\tau^2=1 \text{ et
}\sigma\tau=\tau\sigma.\eqno{(7.1)}$$ On note $\kappa $ le noyau de
l'homomorphisme $\mathcal{N}_{\mathbf{k}/\QQ(\sqrt{2p})}$, d'après
le Lemme précédent on a $$C_{\mathbf{k}, 2}/\kappa \simeq <c>\simeq
\ZZ/4\ZZ .\eqno{(7.2)}$$ Soit $\beta$ l'idéal premier de
$\QQ(\sqrt{2p})$ au-dessus de 2, alors on a
$$2=\beta^2\text{ dans
}\QQ(\sqrt{2p}),\quad\beta=\mathcal{H}^2\text{ dans
}\mathbf{k}\quad\text{et}\quad
\mathcal{N}_{\mathbf{k}/\QQ(\sqrt{2p})}([\mathcal{H}])=[\beta].$$
Par suite la classe de $\beta$ est d'ordre 2 dans $\QQ(\sqrt{2p})$
(car l'unité fondamentale de $\QQ(\sqrt{2p})$ est de norme -1). Ce
qui implique que $[\beta]=c^2$, en particulier
$$\mathcal{N}_{\mathbf{k}/\QQ(\sqrt{2p})}([\mathcal{H}])=c^2.\eqno{(7.3)}$$
Les résultats (7.1), (7.2) et (7.3), permettent de conclure le
Théorème énoncé.
\epr
\section{Preuve du théorème principal}
Reprenons la situation et les notations de la section 5, alors le
groupe $G=Gal(\mathbf{k}_2^{(2)}/\mathbf{k})$ est métacyclique et
non-modulaire. Donnons maintenant une démonstration du Théorème
principal, pour cela nous avons besoin du Lemme suivant :

\bl

Soit $\mathbf{k}=\QQ(\sqrt{2p}, i)$ avec $p$ un nombre premier tel
que $(\frac{2}{p})_4=(\frac{p}{2})_4=-1$. Si on suppose que
$G=\langle a, b\rangle=Gal(\mathbf{k}_2^{(2)}/\mathbf{k})$ tel que
$a^2\equiv b^4\equiv 1\mod G'$, alors $ab^2=b^2a$. \el \bpr

Nous avons vu que la suite des 2-corps de classes de $\mathbf{k}^*$
s'arrête en ${\mathbf{k}^{*}}^{(1)}_2=\mathbf{k}_2^{(2)}$, alors le
groupe $Gal(\mathbf{k}_2^{(2)}/\mathbf{k}^*)=H_{3, 2}=\langle a,
b^2, G'\rangle$ est abélien, ce qui implique que $ab^2=b^2a$. \epr
\bpr[Preuve du théorème principal] Montrons que la classe
$[\mathcal{H}_i]$ de $\mathcal{H}_i$ est d'ordre 2. Comme $\pi_i$ se
ramifie dans $\mathbf{k}/\QQ(i)$ et $p$ se ramifie dans
$\QQ(\sqrt{2p})$, alors il existe un idéal $\mathcal{H}_i$ de
$\mathbf{k}$ tel que $\mathcal{H}_i^2=(\pi_i)$ et un idéal $\beta$
de $\QQ(\sqrt{2p})$ tel que $\beta^2=(p)$. On suppose que
$\mathcal{H}_i=(y)$ pour un certain $y$ de $\mathbf{k}$. On a donc,
en prenant la norme,
$\mathcal{N}_{\mathbf{k}/\QQ(\sqrt{2p})}(\mathcal{H}_i)=(x)$ pour un
certain $x$ de $\QQ(\sqrt{2p})$ et $\beta^2=(x^2)=(p)$. Ce qui est
équivalent à l'existence d'une unité $\epsilon$ de $\QQ(\sqrt{2p})$
telle que $p\epsilon=x^2$, alors $\epsilon$ est égal à
$\pm\epsilon_{2p}$ où $\epsilon_{2p}$ est l'unité fondamentale de
$\QQ(\sqrt{2p})$ ou bien $\pm 1$, ce qui montre que la norme de
$\epsilon_{2p}$ est positive ou bien $p$ est un carré dans
$\QQ(\sqrt{2p})$, alors on a une contradiction, puisque la norme de
$\epsilon_{2p}$ vaut -1 et $\sqrt p\notin\QQ(\sqrt{2p})$.
Il s'ensuit que la classe de $\mathcal{H}_i$ est d'ordre 2.\\
\indent Montrons que $\mathcal{H}_i$ et $\mathcal{H}$ représentent
la même classe dans $\mathbf{k}$. Rappelons que si $L/M$ est une
extension cyclique, on désigne par $\mathrm{Am}(L/M)$ le groupe de classes
ambigues et par $\mathrm{Am}_f(L/M)$ celui des classes fortement ambigues.
Alors on a $$|Am(L/M)|=|Am_f(L/M)|[E_M\cap
\mathcal{N}_{L/M}(L^\times):\mathcal{N}_{L/M}(E_L)],$$ où $E_L$
(resp. $E_M$) est le groupe des unités de $L$ (resp $M$). Dans notre
cas le 2-groupe de classes de $\mathbf{k}$  est de type $(2, 4)$,
donc $|Am(\mathbf{k}/\QQ(i))|=4$. On a aussi $E_{\mathbf{k}}$ est
engendré par $\epsilon_{2p}$ et $i$, or $i$ est norme dans
$\mathbf{k}/\QQ(i)$, par suite $|Am_f(\mathbf{k}/\QQ(i))|=2$,
c'est-à-dire qu'il existe une seule classe de $\mathbf{k}$ d'ordre
2, fortement ambigue. Comme $(\pi_i)=\mathcal{H}_i^2$ et
$(1+i)=\mathcal{H}^2$, alors les classes de $\mathcal{H}_i$ et
$\mathcal{H}$ sont fortement ambigues. Finalement on a
$[\mathcal{H}_1]=[\mathcal{H}_2]=[\mathcal{H}]$.\\
\indent Montrons que $G=Gal(\mathbf{k}_2^{(2)}/\mathbf{k})=\langle
a, b\rangle$ est un groupe métacyclique non-modulaire  où
$a^{2^n}=b^4=1$ et $b^{-1}ab=a^{-1+k2^{s}}$ tels que $1<s<\alpha$,
$k$ un nombre impair. On sait, d'après la Proposition \ref{ta08} que
la classe de $\mathcal{H}$ est d'ordre 2 et capitule dans
$\mathbf{k}^*=\mathbf{K}_{3, 2}$. Or, on a vu, (voir Théorème
\ref{ta10}) qu'il y a exactement deux classes qui capitulent dans
$\mathbf{k}^*$, donc $kerj_{\mathbf{k}\rightarrow
\mathbf{k}^*}=[\mathcal{H}]=\sigma^2$. Par la loi de réciprocité
d'Artin, on trouve que $kerV_{G\rightarrow H_{3, 2}}=b^2G'$. D'après
la Proposition \ref{ta07} et le Lemme précédent, on a
$V_{G\rightarrow H_{3, 2}}(b^2G')=b^{-1}b^2bb^2 H^{'}_{3, 2}=b^4=1$,
ce qui se produit seulement si $G$ est un groupe métacyclique
non-modulaire de type 1 ou 3 (voir \ref{ta06}). Supposons que $G$
est de type 1, donc $b^{-1}ab=a^{-1}$. D'après la proposition
\ref{ta07}, $V_{G\rightarrow H_{3, 2}}(aG')=b^{-1}aba=a^{-1}a=1$, ce
qui implique qu'il y a exactement quatre classes qui capitulent dans
$\mathbf{k}^*$; ce qui n'est pas notre cas. On a donc $G$ est de
type 3 c'est-à-dire $G=\langle a, b\rangle$ où $a^{2^{\alpha}}=1,
b^{4}=1, b^{-1}ab=a^{-1+k2^s}$tels que $1<s<\alpha$, $k$ un nombre
impair et $G$ est d'ordre $2^{\alpha+2}$. Or, on sait, d'après le
Théorème \ref{ta16} que la suite des 2-corps de classes de
$\mathbf{k}^*$ s'arrête en
${\mathbf{k}^{*}}^{(1)}_2=\mathbf{k}_2^{(2)}$, alors
 $Gal({\mathbf{k}^*}_2^{(1)}/\mathbf{k}^*)=Gal(\mathbf{k}_2^{(2)}/\mathbf{k}^*)$.
De plus on a $h(\mathbf{k}^*)=2.h(-p)=2.2^n$, par conséquent l'ordre
de $G$ est égal à $2^{n+2}$ et $\alpha=n$.\\
\indent Montrons que seules la classe de $\mathcal{H}$ est son carré
capitulent dans chacune des trois extensions quadratiques non
ramifiées de $\mathbf{k}$. On a vu que la classe de $\mathcal{H}$
est son carré capitulent dans $\mathbf{k}^*$. Calculons
$V_{G\rightarrow H_{1, 2}}(aG')$, $V_{G\rightarrow H_{1,
2}}(b^2G')$,$V_{G\rightarrow H_{2, 2}}(aG')$ et $V_{G\rightarrow
H_{2, 2}}(b^2G')$. Remarquons d'abord que $H_{1, 2}=\langle b,
G'\rangle=\langle b, a^2\rangle$, donc le groupe des commutateurs
$H^{'}_{1, 2}=\langle a^{-2}b^{-1}a^2b\rangle$. Comme
$b^{-1}ab=a^{-1+k2^{s}}$, alors
$$H^{'}_{1, 2}=\langle a^{-2}a^{-2+k2^{s+1}}\rangle=\langle a^{-2+k2^{s}}\rangle^2=\langle a^{-1}b^{-1}ab\rangle^2=G'^2=\langle a^4\rangle.$$
\indent De la même façon, on trouve que $H^{'}_{2, 2}=\langle
a^{4}\rangle$. Les résultats de la Proposition \ref{ta07} et le
Lemme précédent montrent
 que $$V_{G\rightarrow H_{1, 2}}(b^2G')=a^{-1}b^2ab^2H^{'}_{1, 2}=a^{-1}ab^2b^2H^{'}_{1, 2}=b^4H^{'}_{1, 2}=H^{'}_{1, 2}=1,$$
avec le même raisonnement précédent, on montre les résultats
suivants
$$V_{G\rightarrow H_{1, 2}}(aG')=a^2H^{'}_{1, 2}=V_{G\rightarrow H_{2, 2}}(aG')\quad \text{et}\quad V_{G\rightarrow H_{2, 2}}(b^2G')=H^{'}_{2, 2}=1.$$
\indent Puisque $a^2\notin H^{'}_{1, 2}$ et $H^{'}_{2, 2}$, par la
loi de réciprocité d'Artin, on trouve que seules la classe de
$\mathcal{H}$ est son carré capitulent dans $\mathbf{K_{1,
2}}/\mathbf{k}$ et dans $\mathbf{K_{2,  2}}/\mathbf{k}$.\\
\indent Montrons que les huit classes de $C_{\mathbf{k}, 2}$
capitulent dans les trois extensions abéliennes non ramifiées de
degré $4$ de $\mathbf{k}$. On peut comme précédemment montrer que
$V_{G\rightarrow H_{1, 4}}(aG')=V_{G\rightarrow H_{2,
4}}(aG')=V_{G\rightarrow H_{3, 4}}(aG')=a^{k2^{s+1}}$ et
$V_{G\rightarrow H_{1, 4}}(bG')=V_{G\rightarrow H_{2,
4}}(bG')=V_{G\rightarrow H_{3, 4}}(bG')=b^4$. Soient $F$ le corps de
genres de $\mathbf{k}/\QQ(i)$ et $\mathrm{Am}(\mathbf{k}/\QQ(i))$ le
sous-groupe des classes ambigues dans $\mathbf{k}/\QQ(i))$ de
$C_{\mathbf{k}, 2}$. Comme le nombre de classes de $\QQ(i)$ est égal
à $1$, alors, d'après la théorie du genres, on a
$$\mathrm{Am}(\mathbf{k}/\QQ(i))\simeq Gal(F/\mathbf{k})\simeq
C_{\mathbf{k}, 2}/(C_{\mathbf{k}, 2})^2\simeq \ZZ/2\ZZ\times
\ZZ/2\ZZ,$$ en particulier $\mathrm{Am}(\mathbf{k}/\QQ(i))=\langle\sigma^2,
\tau\rangle$ et $F=\mathbf{K}_{3, 4}$ (car $\mathbf{K}_{3, 4}$ est
la seule extension abélienne non ramifiée de type $(2, 2)$  sur
$\mathbf{k}$). Comme $\langle\sigma^2, \tau\rangle\simeq \langle
b^2G', aG'\rangle$, alors d'après F. Terada, toute les classes
ambigues de $\mathbf{k}$ relativement à $\QQ(i)$ capitulent dans
$\mathbf{K}_{3, 4}$ (voir par exemple \cite{Su91}). Par suite on a
$\langle b^2G', aG'\rangle\subset kerV_{G\rightarrow H_{3, 4}}$. Les
résultats précédents impliquent que $a^{k2^{s+1}}=1$, par conséquent
$2^n$ divise $k2^{s+1}$. Puisque $k$ est un nombre impair et $a$
d'ordre $2^n$ , on trouve que $2^n$ divise $2^{s+1}$ et $n\leq s+1$,
or on a $s\leq n-1$, ce qui prouve que $s=n-1$ et
$a^{k2^{s+1}}=a^{k2^n}=b^4=1$, ainsi
$$V_{G\rightarrow H_{1, 4}}(G/G')=V_{G\rightarrow
H_{2, 4}}(G/G')=V_{G\rightarrow H_{3, 4}}(G/G')=1.$$ Ceci achève la
preuve du théorème principal .\epr

\end{document}